\documentclass[12pt,a4paper]{article}

\usepackage[cp1251]{inputenc}
\usepackage[english]{babel}

\usepackage{amsfonts,amsmath,amsxtra,amsthm,amssymb,latexsym}

\newcommand{\permanent}[1]{\footnote{Permanent address: #1}}
\newcommand{\email}[1]{{\upshape(\texttt{#1})}}

\newcommand{\ack}{\section*{Acknowledgments}}
\theoremstyle{plain}
\newtheorem{theorem}{Theorem}[section]
\newtheorem{corollary}{Corollary}[section]
\newtheorem{lemma}{Lemma}[section]
\newtheorem{problem}{Problem}

\newtheorem{definition}{Definition}[section]
\newtheorem{remark}{Remark}[section]

\newtheorem{proposition}{Proposition}[section]

\numberwithin{equation}{section}

\DeclareMathOperator{\supp}{supp} \DeclareMathOperator{\im}{Im}
\DeclareMathOperator{\Real}{Re} 
\DeclareMathOperator{\sgn}{sgn} \DeclareMathOperator{\dom}{dom}
\DeclareMathOperator{\ran}{ran} 
\DeclareMathOperator{\Span}{span}

\def\R{\mathbb R}
\def\C{\mathbb C}

\def\J{\mathcal{J}}
\def\H{\mathfrak{H}}
\def\HH{\mathcal{H}}
\def\K{\mathcal{K}}
\def\Ainf{\mathcal{A}_{\infty}}
\def\AR{\mathcal{A}_0}
\def\ep{\varepsilon}
\def\m{\mathcal{M}}

\def\Lminm{L_{\min\!-}}
\def\Lminp{L_{\min\!+}}
\def\Lminpm{L_{\min\!\pm}}
\def\Aminm{A_{\min\!-}}
\def\Aminp{A_{\min\!+}}
\def\Aminpm{A_{\min\!\pm}}
\def\Azm{A_{0-}}
\def\Azp{A_{0+}}
\def\A0pm{A_{0\pm}}
\def\Lzm{L_{0-}}
\def\Lzp{L_{0+}}
\def\len{l}
\def\L{\mathcal{L}}

\title{Indefinite Sturm-Liouville  operators with the singular critical
point zero}

\author{Illya~M.~Karabash\permanent{Department of Partial Differential Equations, Institute
of Applied Mathematics and Mechanics of NAS of Ukraine, R.
Luxemburg str., 74, Donetsk 83114, UKRAINE
\email{karabashi@mail.ru}}\ \ and
Aleksey~S.~Kostenko\permanent{Department of Nonlinear Analysis,
Institute of Applied Mathematics and Mechanics of NAS of Ukraine,
R. Luxemburg str., 74, Donetsk 83114, UKRAINE
\email{duzer80@mail.ru}}}

\date{}
\begin{document}

\maketitle

\begin{abstract}
We present a new necessary condition for similarity of indefinite
Sturm-Liouville operators to self-adjoint operators. This
condition is formulated in terms of Weyl-Titchmarsh $m$-functions.
Also we obtain necessary conditions for regularity of the critical
points $0$ and $\infty$ of $J$-nonnegative Sturm-Liouville
operators. Using this result, we construct several examples of
operators with the singular critical point zero. In particular, it
is shown that $0$ is a singular critical point of the operator
$-\frac{(\sgn x)}{(3|x|+1)^{-4/3}} \, \frac{d^2}{dx^2}$\ \ acting
in the Hilbert space $L^2(\R, (3|x|+1)^{-4/3}dx)$ and therefore
this operator is not similar to a self-adjoint one. Also we
construct a J-nonnegative Sturm-Liouville operator of type $(\sgn
x)(-d^2/dx^2+q(x))$ with the same properties.
\end{abstract}

\section{ Introduction}\label{intro}

In this paper, we are interested in Sturm-Liouville equations
\begin{equation}\label{1}
-y''(x)+q(x)y(x)=\lambda\ r(x)y(x), \qquad x\in \R,
\end{equation}
with an indefinite weight $r$. More specifically, we study the
spectral properties of the associated non-self-adjoint operator
\begin{equation}\label{3}
A := \frac 1{r} \left( -\frac{d^2}{dx^2} +q \right)
\end{equation}
acting in the weighted Hilbert space $L^2(\R, \ |r(x)|dx)$ (an
explicit definition of the operator $A$ is given in the next
section). Here the weight $r$ and the potential $q$ are real and
locally Lebesgue integrable functions on $\R$ \ ($q, r \in
L_{loc}^1(\R)$), and $x r(x)> 0$ for all $x\in \R\setminus\{0\}$.
Thus $r$ changes sign at $0$.

The spectral problem
\begin{equation}\label{e eq|r|}
-y''(x)+q(x)y(x)=\lambda |r(x)| y(x), \qquad x\in \R,
\end{equation}
with the positive weight $|r|$ is usually treated in the context
of the Hilbert space $L^2 (\R, |r(x)|dx)$ with the scalar product
$(f,g)= \int_R f\overline{g} |r| dx$. Under the assumption that
\eqref{e eq|r|} is in the limit point case at $-\infty$ and
$+\infty$, the operator $L$ associated with \eqref{e eq|r|} is
self-adjoint in $L^2 (\R, |r(x)|dx)$ and the operator $A$
associated with \eqref{1} is \emph{J-self-adjoint}. The letter
means that $A$ is self-adjoint with respect to the indefinite
inner product
\[
[f,g]:=(Jf,g)= \int_R f\overline{g}\, r \, dx ,
\]
where the operator $J$ is defined by $J:f(x) \mapsto (\sgn x)
f(x). $
% is a fundamental
%symmetry in the Krein space $L^2 (\R, r(x)dx)$ (see \cite{AzJ89},
%\cite{Lan82}).
Obviously, the operators $A$ and $L$ are connected by the equality
$A=JL$. Notice also that the operator $A$ is non-self-adjoint (in
the Hilbert space $L^2 (\R, |r(x)|dx)$).

The main object of the present paper is the \emph{similarity} of
the operator $A$ to a self-adjoint operator. Let us recall that
two closed operators $T_1$ and $T_2$ in a Hilbert space $\H$ are
called \emph{similar} if there exist a bounded operator $S$ with
the bounded inverse $S^{-1}$ in $\H$ such that $S\dom(T_1) =
\dom(T_2)$ and $T_2 = S T_1 S^{-1}$.

The similarity of the corresponding J-self-adjoint operators to a
self-adjoint operator is essential for the solution of
forward-backward boundary value problems, which arise in certain
physical models, particularly in transport and scattering theory
(see \cite{B85,KLH82,GvdMP87,KarKr06,KarGAMM06}), and in the theory of random
processes (see \cite{Pag74} and references therein).

If the operator $L$ is nonnegative, $L \geq 0$, one can study the
similarity problem for the operator $A$ in the context of the
spectral theory of \emph{J-nonnegative operators} \cite{Lan82}
(necessary notions and facts are contained in Section
\ref{indef}). If, in addition, the resolvent set of the operator
$A$ is nonempty, $\rho(A)\neq\emptyset$, then the operator $A$
possesses the following properties:
\begin{description}
\item[(i)] the spectrum of $A$ is real, $\sigma(A) \subset \R$;
\item[(ii)] if $\lambda \neq 0$ is an eigenvalue of $A$, then it
is semisimple (i.e., $\ker (A-\lambda) = \ker (A-\lambda)^2$);
\item[(iii)] if $0 $ is an eigenvalue of $A$, then its Riesz index
$\leq 2$, i.e., $\ker A^2 = \ker A^3$ (generally, $0$ may be a
nonsemisimple eigenvalue).
\end{description}
 Moreover,
$A$ admits a \emph{spectral function} $E_A (\Delta)$. The
properties of $E_A (\Delta)$ are similar to the properties of a
spectral function of a self-adjoint operator. The main difference
is the occurrence of \emph{critical points}. Significantly
different behavior of the spectral function $E_A (\Delta)$ occurs
at \emph{singular critical point} in any neighborhood of which the
spectral function is unbounded. The critical points, which are not
singular, are called \emph{regular}. It should be stressed that
only $0$ and $\infty$ may be critical points for J-nonnegative
operators.
%Besides, $A$ may have a nonsemisimple eigenvalue at $0$.
Under the additional assumption $\ker A = \ker A^2$,
%(i.e., either $0$ is a semisimple eigenvalue of $T$ or $0 \in \rho(T)$)
the following assertions  are equivalent: 
%(see for example\cite{Lan82}):
\begin{description}
\item[(i)] $A$ is similar to a self-adjoint operator in $\H$;
\item[(ii)]  $0$ and $\infty$ are regular critical points of $A$.
\end{description}

%If the spectrum $\sigma(A)$ is real and discrete, the similarity
%of $A$ to a self-adjoint operator is equivalent to the Riesz basis
%property of eigenvectors.
In \cite{B85}, Beals showed that the eigenfunctions of regular
Sturm-Liouville problems 
of the type \eqref{1} 
form a Riesz basis
if $r(x)$ behaves like a power of $x$ at $0$. Improved versions of
Beals' condition have been provided in
\cite{CurLan,F96,Par03,Pyat89,Sh93,Vol96}. In \cite{CurLan,F96},
singular differential operators with indefinite weights have been
considered and the regularity of the critical point $\infty$ was
proven for a wide class of weight functions. The existence of
indefinite Sturm-Liouville operators with the singular critical
point $\infty$ was established in \cite{Vol96}, and corresponding
examples were constructed in \cite{AbPyat97,F98,Par03}. The
question of nonsingularity of $0$ is much harder. It was shown in
\cite{CN95,FSh1,FN98,KarMFAT00,Kar_Mal,KM06,Kos05,Kos06,Kos_OT_06} that
$0$ is a regular critical point for several model classes of
differential operators. In \cite{KM06} several  necessary
similarity
conditions in terms of Weyl functions
were obtained also. The following problem naturally arises in
this context:
\begin{problem}
%\textbf{Problem.}
Whether there are any J-nonnegative Sturm-Liouville operators $A$
with the singular critical point $0$.
\end{problem}

 It will be shown in the present paper that those operators
do exist.

The paper is organized as follows. In Section \ref{prelim} we
summarize necessary definitions and statements from the spectral
theory of Sturm-Liouville operators and from the spectral theory
of $J$-nonnegative operators.

The main results of the paper are contained in Section \ref{ness}.
The central result is a necessary condition for the operator $A$
to be similar to self-adjoint one (Theorem \ref{necessity}).
Further, we obtain necessary conditions for regularity of the
critical points $0$ and $\infty$ of $J$-nonnegative
Sturm-Liouville operators (Theorem \ref{thIV.1}). These conditions
are formulated in terms of the Weyl-Titchmarsh $m$-functions.
Proofs of these results are contained in Section \ref{s proof}.

In Section \ref{s ex}, we show that $0$ is a singular critical
point of the operator associated with the differential expression
\begin{gather*}
 -\frac{\sgn x}{(3|x|+1)^{-4/3}} \frac{d^2 }{dx^2}
\end{gather*}
in the Hilbert space $L^2(\R, (3|x|+1)^{-4/3})$. Moreover, we
construct a J-nonnegative operator $A_0=(\sgn x) (-d^2/dx^2 +q_0)$
with the same property in Subsection \ref{sub_6_2}.

In \cite{FSh1}, Faddeev and Shterenberg proved the following:
\emph{
if the operator $A=(\sgn x)(-d^2/dx^2+q)$ (acting in $L^2(\R)$) has
a real spectrum and
\begin{gather}
\int_\R |q(x)| (1+|x|^2)\ dx < \infty \label{e fmoment},
\end{gather}
then the operator $A$ is similar to a self-adjoint one}. In
Subsection \ref{sub_6_1}, we show that \eqref{e fmoment} cannot be
changed to the condition $\int_\R |q(x)| (1+|x|^\alpha ) dx <
\infty$ with $\alpha <1$.

It should be noted that the necessary conditions obtained in
\cite{KM06} ignore the singular part of the operator $A$.
Moreover, these conditions are fulfilled for examples constructed
in Sections 5 and 6. But connections between \cite[Corollaries
5.4--5.6]{KM06} and Theorem \ref{necessity} are not clear (see
also Remark \ref{rIII.1}).

\textbf{Notation:} $\mathfrak{H}, \mathcal{H}$ denote separable
Hilbert spaces. The scalar product and the norm in the Hilbert
space $\mathfrak{H}$ are denoted by $(\cdot,\cdot)_{\mathfrak{H}}$
and $\|\cdot \|_{\mathfrak{H}}$, respectively. The set of all
bounded linear operators from $\mathfrak{H}$ to $\mathcal{H}$ is
denoted by $[\mathfrak{H}, \mathcal{H}]$ or $[\mathfrak{H}]$ if
$\mathfrak{H}= \mathcal{H}$. $ \Span \{ f_1,f_2, \dots, f_N \}$
denotes the closed linear hull of vectors $f_1$, $f_2$, \dots,
$f_N$.
%$\mathcal{C}(\mathfrak{H})$ stands
%for the set of closed densely defined operators in $\mathfrak{H}$.
Let $T$ be a linear operator in a Hilbert space $\mathfrak{H}$. In
what follows $\dom (T)$, $\ker (T)$, $\ran (T)$ are the domain,
kernel, range of $T$, respectively. We denote by $\sigma(T)$,
$\rho (T)$ the spectrum and the resolvent set of $T$; $\sigma_p
(T)$ stands for the set of eigenvalues of $T$. $ R_T \left(
\lambda \right):=\left( T-\lambda I\right)^{-1} $, $\lambda \in
\rho(T)$, is the resolvent of $T$.

We set $\C_\pm := \{ \lambda \in \C \ : \ \pm \im \lambda
>0\}$. By
$\chi_\pm(t):=\chi_{\R_\pm }(t)$ we denote the characteristic
function of $\R_\pm $, where $\R_+:=[0, +\infty)$, $
\R_-:=(-\infty, 0]$. We write $f\in L^1_{loc}(\R) (\in
AC_{loc}(\R))$ if the function $f$ is Lebesgue integrable
(absolutely continuous) on every bounded interval in $\R$.
%We denote by $L^2(\R, |r(x)|dx)$ the
%Hilbert function space of equivalence classes of Lebesgue
%measurable functions $f$ such that $\int_{\R}
%|f(x)|^2|r(x)|dx<\infty$; the inner-product of $ f, g\in L^2(\R,
%|r(x)|dx)$ is defined by $(f, g):=\int_{\R}
%f(x)\overline{g(x)}|r(x)|dx$.

 \section{Preliminaries}\label{prelim}

\subsection{Differential operators.}\label{ss DifOp}

Consider the differential expressions
\begin{equation}\label{II_1_01}
\ell [y] = \frac 1{|r|} \left( -y'' +qy \right) \qquad \text{and}
\qquad a [y] =  \frac 1{r} \left( -y'' +qy \right).
\end{equation}
Under the assumptions
%the weight function $r$ and the potential
%$q$ are locally Lebesgue integrable on $\R$,
$q, r \in L_{loc}^1(\R)$, these differential expressions are
regular on all compact intervals $[a, b]\subset \R$, but singular
on $(-\infty, +\infty)$.

Let $\mathfrak{D}$ be the set of all functions $f \in L^2 (\R,
|r(x)|dx)$ such that $f$ and its first derivative $f'$ are locally
absolutely continuous and  $\ell [f] \in L^2 (\R, |r(x)|dx)$,
\begin{equation}\label{II_1_02}
\mathfrak{D}:=\{f \in L^2 (\R, |r(x)|dx):\ f, f'\in AC_{loc}(\R),
\ \ell [f] \in L^2 (\R, |r(x)|dx)\}.
\end{equation}
The set $\mathfrak{D}$ is the maximal linear manifold in $L^2 (\R,
|r(x)|dx)$ on which the differential expressions $\ell
[\cdot]$ and $a[\cdot]$ have a natural meaning. %and $\ell [y] \in L^2 (\R)$.
%consists of $y \in L^2 (\R, |r|)$ which
On $\mathfrak{D}$ let us define the operators $L$ and $A$ as
follows:
\begin{gather}
\dom(L) = \dom(A) = \mathfrak{D}, \notag \\
L f = \ell [f] , \qquad A f =  a[f]\qquad \text{for} \qquad f\in
\mathfrak{D}.\label{II_1_03}
\end{gather}
The operators $A$ and $L$ are called \emph{the operators
associated with equations} \eqref{1} and \eqref{e eq|r|},
respectively. It is well known (see \cite{AG}) that $A$ and $L$
are closed differential operators in $L^2 (\R, |r(x)|dx)$.

  Moreover, the adjoint operator $L^*$ of $L$ is a
closed symmetric operator with deficiency indices $(n, n), \ 0\leq
n\leq 2$. In what follows we always assume that the differential
expression $\ell [\cdot]$ is \emph{in the limit point case at
$+\infty$ and $-\infty$}. In other words, we assume that $n=0$,
i.e., the operator $L$ is self-adjoint in the Hilbert space $L^2
(\R, |r(x)|dx)$.
%Hence $A$ is a J-self-adjoint
%operator. The operator $A$ is called \emph{the  operator
%associated with equation} \eqref{1}.

It is clear that $A=JL$, where the operator $J$ is defined by
\begin{gather*}
(Jf)(x) = (\sgn x) f (x) , \qquad f \in L^2 (\R, |r(x)|dx).
\end{gather*}
Obviously, $J^*=J^{-1}=J\ $ \ in $L^2 (\R, |r(x)|dx)$, and
$A^*=LJ$. Thus the adjoint operator $A^*$ of $A$ is defined by the
differential expression $a [\cdot]$ on the domain
\begin{equation}
 \dom (A^*) = J^{-1} \mathfrak{D} = J \mathfrak{D} = \{ f \in L^2 (\R, |r(x)|dx ) \ : \ Jf \in
\mathfrak{D} \}.
\end{equation}
Since $\dom (A^*) \neq \dom (A)$, we have $A \neq A^*$, i. e., the
operator $A$ is non-self-adjoint in $L^2 (\R, |r(x)|dx )$.

Let us determine the following set
\begin{equation}
\mathfrak{D}_0:=\dom (A)\cap\dom (A^*)=\mathfrak{D}\cap
J\mathfrak{D}=\{ f \in \mathfrak{D}\ : f(0)=f'(0) = 0 \}.
\end{equation}
It is obvious that the following restrictions of the operators $L$
and $A$
\begin{equation} \label{e Amin=}
L_{\min}:=L\upharpoonright \mathfrak{D}_0,\qquad
A_{\min}:=A\upharpoonright \mathfrak{D}_0.
\end{equation}
are closed densely defined symmetric operators in $L^2 (\R,
|r(x)|dx )$ with equal deficiency indices $n_\pm(L_{
\min})=n_\pm(A_{ \min})=2$. Let $P_\pm$ denote the orthogonal
projections in $L^2 (\R, |r(x)|dx )$ onto $L^2 (\R_\pm, |r(x)|dx
)$. One can represent the operator $L_{\min}$ in the
following form
\begin{equation}\label{II_I_7}
 L_{ \min}= \Lminm \oplus \Lminp,
\end{equation}
where
\begin{equation}\label{II_I_8}
\Lminpm:=L\upharpoonright L^2(\R_\pm, |r|dx) \ = \ L\upharpoonright
\mathfrak{D}_0^\pm;\qquad
\mathfrak{D}_0^{\pm}=P_\pm \mathfrak{D}_0.
\end{equation}

The operators $\Lminp$ and $\Lminm$ are called \emph{the
minimal operators associated with the differential expression
$\ell [\cdot]$} on $\R_+$ and $\R_-$, respectively (see
\cite{AG}). Notice that $\Lminpm$ is a symmetric operator in
the Hilbert space $L^2 (\R_\pm, |r(x)|dx )$. The deficiency
indices of $\Lminpm$ are equal to (1,1). In what follows
$\mathfrak{D}_0^*$ stands for the domain of the adjoint operator
$L_{\min}^*$ of $L_{\min}$.

It is not hard to see that
\begin{equation}\label{II_I_9}
A_{ \min}= \Aminm \oplus \Aminp, \qquad \text{where} \quad \Aminpm
:=\pm \Lminpm .
\end{equation}
Thus the operator $A$ is a non-self-adjoint extension of the
operator $A_{\min}$ and
\begin{multline} \label{e domA}
\mathfrak{D} = \dom (A) = \left\{ f \in \dom(\Aminp^*)
\oplus \dom(\Aminm^*):
%\right. \\ \left.
\ f(+0) = f(-0),\ \ f'(+0) = f' (-0) \right\} 
\end{multline}
(note that 
$\dom(\Aminpm^*)=P_\pm \dom(A_{\min}^*) =P_\pm \mathfrak{D}_0^* $~).

\subsection{Weyl-Titchmarsh m-coefficients.}\label{sII_2}

Let $c(x, \lambda)$ and $s(x, \lambda)$ denote the linearly
independent solutions of equation \eqref{e eq|r|} satisfying the
following initial conditions at zero
\[
c(0, \lambda)= s'(0, \lambda)=1; \qquad c^{\prime}(0, \lambda)=
s(0, \lambda)=0.
\]
Since equation \eqref{e eq|r|} is limit-point at $+\infty$, the
Weyl-Titchmarsh theory states (see, for example, \cite{LevSar70})
that there exists a unique holomorphic function
$m_+(\cdot):\C\setminus\R\to\C$, such that the function $s(x,
\lambda) - m_+ (\lambda) c(x, \lambda) $ belongs to $L^2(\R_+,
|r(x)|dx)$. Similarly, the limit point case at $-\infty$ yields
the fact that there exists a unique holomorphic function
$m_-(\cdot):\C\setminus\R\to\C$, such that $ s(x, \lambda) + m_-
(\lambda) c(x, \lambda)\in L^2(\R_-, |r(x)|dx)$. (Note that for
$\lambda \in \C \setminus \R$ the functions $c(x,
\lambda)\chi_\pm(x)$ and $s(x, \lambda)\chi_\pm(x)$ do not belong
to $L^2(\R, |r(x)|dx)$).

The functions $m_+$ and $m_-$ are called \emph{the Weyl-Titchmarsh
m-coefficients for} \eqref{e eq|r|} on $\R_+$ and on $\R_-$,
respectively. We put
\begin{equation}\label{e def psi}
M_\pm (\lambda) := \pm m_\pm (\pm \lambda);\qquad \psi_\pm(x,
\lambda)=(s(x, \pm\lambda)-M_\pm(\lambda)c(x,
\pm\lambda))\chi_\pm(x).
\end{equation}
By the definition of $m_\pm$, the functions
$\psi_+ (\cdot , \lambda)$ and $\psi_- (\cdot , \lambda)$ belong
to $L^2 (\R, |r(x)| dx)$ for all $\lambda \in \C \setminus \R$. Besides,
$$
a [\psi_\pm(x, \lambda)] = \lambda \psi_\pm(x, \lambda) .
$$
 The function
$M_+ (\cdot)$ ($M_- (\cdot)$) is said to be \emph{the
Weyl-Titchmarsh m-coefficient for} equation \eqref{1} on $\R_+$
(on $\R_-$).

\begin{definition}[e.g. \cite{KK1}]
The class $(R)$ consists of all holomorphic functions $G :\C_+
\cup \C_- \rightarrow \C$ such that
\[
G (\overline{\lambda})= \overline{G(\lambda)}, \qquad \text{and}
\qquad \im\lambda \cdot \im G(\lambda) \geq 0  \qquad \text{for}
\quad \lambda \in \C_+ \cup \C_-.
\]
\end{definition}

It is well known that
\begin{equation}\label{217}
\int_0^{+\infty} | \psi_+(x, \lambda) |^2  r(x) dx =  \frac {\im
M_+ (\lambda)}{\im \lambda} , \qquad \int_{-\infty}^0 | \psi_-(x,
\lambda) |^2  |r(x)| dx =  \frac {\im M_- (\lambda)}{\im \lambda}
,
\end{equation}
for all $\lambda \in \C\setminus\R$ (see e.g. \cite{LevSar70}). These formulae imply that the
functions $M_+$ and $M_-$ (as well as $m_+$ and $m_-$) belong to
the class $(R)$.
Moreover (see \cite{krein53b}, \cite{KK2}, and also \cite[\S II.5,
Theorem 5.2]{LevSar70} for the case $|r|\equiv 1$), the functions
$M_+$ and $M_-$ admit the following integral representation
\begin{equation}\label{216}
M_\pm (\lambda) = \int_{-\infty}^{+\infty}\frac{ d\tau_\pm
(s)}{s-\lambda}, \qquad \lambda \in \C\setminus\R.
\end{equation}
Here  %$\alpha_\pm \in \R$, and
$\tau_\pm:\R\to\R$ are nondecreasing functions on $\R$ with the
following properties
\begin{gather*} %\label{21herglotz}
\int_{-\infty}^{+\infty}\frac{d\tau_\pm (s)} {1+|s|}<+\infty,
\qquad \tau_\pm (0)= 0, \qquad \tau_\pm(s) = \tau_\pm(s - 0) .
\end{gather*}
 Notice that the functions $\tau_+$ and $\tau_-$ are uniquely determined by
the Stieltjes inversion formula (see \cite{KK1}),
\begin{gather} \label{e St}
\lim_{\varepsilon\downarrow 0}\frac{1}{\pi}\int_0^{s}\im M_\pm
(t+i\varepsilon)dt = \frac{\tau_\pm (s+0)+\tau_\pm (s-0)}{2} \ .
\end{gather}
The function $\tau_\pm $ is called a \emph{spectral function of
the boundary value problem}
\begin{gather}\label{s_f}
-y''(x) +q(x) y(x) = \lambda r(x) y(x) , \qquad x \in \R_\pm;
\qquad y'(\pm 0) = 0.
\end{gather}
In other words, the self-adjoint operator
\begin{equation}\label{A0_pm}
\A0pm:=\Aminpm^* \upharpoonright
\{y\in\dom(\Aminpm^*): y'(\pm 0)=0\},
\end{equation}
associated with \eqref{s_f} is unitary equivalent to the
multiplication by the independent variable in the Hilbert spaces
$L^2(\R, d\tau_\pm(x))$. This fact obviously implies
\begin{equation} \label{e sA=supp}
\sigma(\A0pm)=\supp(d\tau_\pm).
\end{equation}
Here $\supp d\tau$ denotes the
\emph{topological support} of a Borel measure $d\tau$ on $\R$,
i.e., $\supp d\tau$ is the smallest closed set $\Omega\subset\R$
such that $d\tau (\R \setminus \Omega) = 0$.

\begin{remark}\label{remII_I_1}
 It is well known that the functions
$M^\infty_\pm:=-1/M_\pm$ belong to the class $(R)$ also. Besides,
they admit the following integral representation
\[
M^\infty_\pm(\lambda)=C_\pm +
\int_{-\infty}^{+\infty}\left(\frac{1}{s-\lambda}-\frac{s}{1+s^2}\right)d\tau^\infty_\pm(s),
\qquad \lambda\in\C\setminus\R,
\]
where $\tau^\infty_+:\R\to\R$ and
$\tau^\infty_-:\R\to\R$ are nondecreasing functions with the following properties
\begin{gather*} %\label{21herglotz}
\int_{-\infty}^{+\infty}\frac{d\tau^\infty_\pm (s)}
{1+|s|^2}<+\infty, \qquad \tau^\infty_\pm (0)= 0, \qquad
\tau^\infty_\pm(s) = \tau^\infty_\pm(s - 0) .
\end{gather*}
The functions $\tau^\infty_+$ and $\tau^\infty_-$ are called the
spectral functions of the boundary value problems
\begin{gather}\label{s_f_D}
-y''(x) +q(x) y(x) = \lambda r(x) y(x) , \qquad x \in \R_\pm;
\qquad y(\pm 0) = 0.
\end{gather}
Sometimes we will say that the functions $M_\pm(\cdot)$ and $M^\infty_\pm(\cdot)$ are
the Weyl-Titchmarsh $m$-coefficients for the boundary
value problems \eqref{s_f} and \eqref{s_f_D}, respectively.
\end{remark}

In the proofs of Theorems \ref{thVI_1} and \ref{thII_7_2} the
following description of $\sigma_p(A)$ will be used.

\begin{proposition}[\cite{Kar_thesis,KarKr06}] \label{p RindEig}
Let $\lambda \in \R$. Assume that
\begin{equation} \label{e tau=tau}
\tau_\pm (\lambda+0) = \tau_\pm (\lambda-0) \quad \text{and}\quad
\int_{\R } \frac 1{|s-\lambda |^2} \ d \tau_\pm (s) < \infty .
\end{equation}
Then $\lambda \in \sigma_p (A)$ if and only if
\begin{gather} \label{e fr1=fr1}
\int_\R \frac 1{s-\lambda} d\tau_+ (t) =
\int_\R \frac 1{s-\lambda} d\tau_- (t) .
\end{gather}
If \eqref{e fr1=fr1} holds true, then:
\begin{description}
\item[(i)] the geometric multiplicity of $\lambda$
is 1;
\item[(ii)] $\lambda$ is a simple eigenvalue if and only if at least one
the following conditions fails
\begin{gather}
\int_{\R } \frac 1{|s-\lambda |^{4}}
\ d \tau_- (t) < \infty ,
\qquad
\int_{\R } \frac 1{|s-\lambda |^{4}} \ d \tau_+ (s) < \infty ,
\label{e frj<inf r}
\\
\int_{\R } \frac 1{(s-\lambda )^{2}} \ d\tau_- (t) =
\int_{\R } \frac 1{(s-\lambda )^{2}} \ d\tau_+ (t)
\label{e frj=frj r} .
\end{gather}
\end{description}
\end{proposition}

\begin{remark}
A complete description of eigenvalues of the indefinite
Sturm-Liouville operator $A$ as well as their geometric and
algebraic multiplicities were obtained in \cite{Kar_thesis}. These
results were published in \cite{KarKr06} and used in \cite{KM06}. Note also
that Proposition \ref{p RindEig} is a particular case of
\cite[Theorem 4.2(3)]{KM06}.
\end{remark}

\subsection{Spectral functions of $J$-nonnegative operators.}\label{indef}

In this subsection basic facts from the spectral theory of
J-nonnegative operators are collected (the reader can find more
details in \cite{AzJ89,Lan82}).

Consider a Hilbert space $\H$ with a scalar product $(\cdot, \cdot
)_\H$. A Hermitian sesquilinear form $ \langle \cdot \, ,\cdot
\rangle $ on a Hilbert space $\H $ is said to be \emph{positive
definite } if $\langle f ,f \rangle >0$ for $f \in
\H\setminus\{0\} $, % $f\neq 0$,
and it said to be \emph{indefinite} if there exist elements $f,g
\in \H$ such that $\langle f,f \rangle <0$ and $\langle g,g
\rangle >0$.

Suppose that $\H = \H_+  \oplus \H_- $, where $\H_+$ and $\H_-$
are closed subspaces of $\H$. Denote by $P_\pm$ the orthogonal
projections from $\H$ onto $\H_\pm$. Let $ \J=P_+ - P_-$ and $[
\cdot, \cdot ]_\K := (\J \cdot ,\cdot )_\H $. Then the pair $\K =
(\H , [\cdot, \cdot]_\K)$ is called a \emph{Krein space} (see
\cite{Lan82,AzJ89} for the original definition). The form $[
\cdot, \cdot ]$ is called \emph{an inner product} in the Krein
space $\K$ and the operator $\J$ is called \emph{a fundamental
symmetry in the Krein space} $\K$. Evidently, the form
$[\cdot,\cdot]$ is indefinite on $\H$ if and only if $\H_+ \neq \{
0\}$ and $\H_- \neq \{ 0\}$.

Let $T$ be a closed densely defined operator in $\H$. The
\emph{J-adjoint operator of} $T$ is defined by the relation
\[
 \ [Tf,g] =
[f,T^{[*]}g] \ , \qquad  f \in \dom (T),
\]
on the set of all $g \in \H$ such that the mapping $f \mapsto
[Tf,g]$ is a continuous linear functional on $\dom (T)$. The
operator $T$ is called \emph{J-self-adjoint} if $ T = T^{[*]} $.
It is easy to see that $T^{[*]} := \J T^* \J$ and the operator $T$
is J-self-adjoint if and only if $\J T$ is self-adjoint. Note that
$\J=\J^*=\J^{-1}=\J^{[*]} $. A closed operator $T$ is called
\emph{J-nonnegative} if $\ [Tf,f] \geq 0 \ $ for $\ f \in \dom (T)
$ (it is equivalent to $\J T \geq 0$).

Let $\mathfrak{S}$ be the semiring consisting of all bounded
intervals with endpoints different from $0$ and $\pm \infty$ and
their complements in $\overline \R := \R \cup{ \infty }$.

\begin{theorem}[\cite{Lan82}] \label{thII.4.1}
Let $T$ be a J-nonnegative J-self-adjoint operator in
$\mathfrak{H}$ with a nonempty resolvent set $\rho (T) \neq
\emptyset$. Then:
\begin{description}
\item[(i)] The spectrum of $T$ is real, $\sigma(T) \subset \R$.
\item[(ii)] There exist a mapping $\Delta \rightarrow E(\Delta)$
from $\mathfrak{S}$ into the set of bounded linear operators in
$\mathfrak{H}$ with the following properties ($\Delta ,\Delta' \in
\mathfrak{S}$):
\begin{description}
\item[(E1)] $E(\Delta \cap \Delta' ) = E(\Delta) E(\Delta' ) $,
\quad $E(\emptyset ) = 0, \quad E(\overline{\R}) = I $, \quad
$E(\Delta) = E(\Delta)^{[*]}$; \item[(E2)] $E(\Delta \cup \Delta'
) = E(\Delta) + E(\Delta' )$ \quad if \quad $\Delta \cap \Delta' =
\emptyset$; \item[(E3)] the form $ \pm [ \cdot,\cdot] $ is
positive definite on $E (\Delta) \mathfrak{H} $, if $\Delta
\subset \R_\pm$; \item[(E4)] $E(\Delta) $ is in the double
commutant of the resolvent of $T$ and $\sigma (T \upharpoonright
E(\Delta) \mathfrak{H}) \subset \overline{\Delta}$; \item[(E5)] if
$\Delta$ is bounded, then $E(\Delta) \mathfrak{H} \subset \dom(T)$
and $T\upharpoonright E (\Delta) \mathfrak{H} $ is a bounded
operator.
\end{description}
\end{description}
\end{theorem}

According to \cite[Proposition II.4.2]{Lan82}, a number $s \in
\{0,\infty\}$ is called \emph{a critical point} of $T$, if the
form $[\cdot,\cdot ]$ is indefinite on $E(\Delta)\mathfrak{H}$ for
each $\Delta \in \mathfrak{S}$ such that $s\in \Delta$. The set of
critical points is denoted by $c(T)$.

If $\alpha \not \in c (T)$, then for arbitrary $\lambda_0, \lambda_1
\in \R \setminus c(T)$, $\lambda_0 < \alpha$, $\lambda_1>\alpha$,
the limits
\begin{equation} \label{KARA_e Ell}
\lim_{\lambda \uparrow \alpha} E([\lambda_0,\lambda]) , \qquad
\lim_{\lambda \downarrow \alpha} E([\lambda,\lambda_1])
\end{equation}
exist in the strong operator topology. % (s.o.t.).
If $\alpha \in c(T)$ and the limits \eqref{KARA_e Ell} do still
exist, then $\alpha$ is called \emph{regular critical point} of
$T$, otherwise $\alpha$ is called \emph{singular}. Here we agree
that, if $\alpha = \infty$, then $\lambda_1 > \alpha$ ($\lambda
\downarrow \alpha$) means $\lambda_1 > -\infty$ ($\lambda
\downarrow -\infty$, respectively).
%\begin{theorem}[\cite{Lan82}]\label{thII.4.2} The homomorphism
%$\Delta \rightarrow E(\Delta)$ can be extended to the semiring
%generated by those intervals whose endpoints are not singular
%critical points of $T$. For this extension properties
%\textbf{(E1)} - \textbf{(E5)} of Theorem \ref{thII.4.1} are
%preserved.
%\end{theorem}

The following proposition is well known.
%and follows from Theorems \ref{thII.4.1}--\ref{thII.4.2}.
\begin{proposition} \label{p cr=sim}
Let T be a J-nonnegative J-self-adjoint operator in the Hilbert
space $\mathfrak{H}$. Assume that  $\rho(T) \neq \emptyset$ and
$\ker T =\ker T^2$ (i.e., $0$ is either a semisimple eigenvalue or
a regular point of T). Then two following assertions are
equivalent:
\begin{description}
\item[(i)] T is similar to a self-adjoint operator. \item[(ii)]
$0$ and $\infty$ are not singular critical points of $T$.
\end{description}
\end{proposition}

\section{Main results}\label{ness}

Let $A$, $L$, and $J$ be the operators defined in Section \ref{ss
DifOp}. Let $M_+$ and $M_-$ be the Weyl-Titchmarsh
$m$-coefficients for \eqref{1}. Throughout this section we always
assume that the operator $A$ is $J$-self-adjoint (this is
equivalent to the self-adjointness of the operator $L$).

In this Section we formulate the main results. All proofs are
contained in the next section.

\subsection{Necessary similarity
condition.}\label{sIII_1}

We start with the following
\begin{proposition}[cf. \cite{KM06}] \label{p NRealSp}
%Let $\lambda \not \in \R$. Then:
\begin{description}
\item[(i)] If $\lambda \not \in \R$, then $\lambda\in \sigma (A)$
if and only if $ \ M_+ (\lambda) = M_- (\lambda ) $.

\item[(ii)] If the operator $L$ is semibounded, then $\rho (A)
\neq \emptyset\ $.
\end{description}
\end{proposition}

\begin{remark}
Proposition \ref{p NRealSp} was obtained in
\cite{KM06} for the case $r(x) = \sgn x$. 
The proof of statement (i) remains the same. 
However, for the operator $A=(\sgn x)(-d^2/dx^2+q)$ with $q \in L^1_{loc} (\R)$,
statement (i) holds true without the assumption $L \geq \eta > -\infty$ 
(see \cite[Proposition 2.5 (v)]{KM06}).
\end{remark}

Let us note that if $A$ is similar to a self-adjoint operator,
then $\sigma(A)\subset\R$ and hence, by Proposition \ref{p
NRealSp}, $M_+ (\lambda) \neq M_- (\lambda)$ for $\lambda \not \in
\R$.

 The central result of the paper is the following theorem.
\begin{theorem}\label{necessity}
If $A$ is similar to a self-adjoint operator, then the functions
\begin{equation}\label{302}
\frac{\im M_+ (\lambda )} {M_+(\lambda )-M_-(\lambda)}\quad
\text{and} \quad \frac{\im
M_-(\lambda)}{M_+(\lambda)-M_-(\lambda)}
\end{equation}
are well-defined and bounded on  $\C\setminus\R$.
\end{theorem}
\begin{corollary}\label{cor_III_1}
Let $A$ be a $J$-self-adjoint operator associated with \eqref{1}.
Assume also that the weight $r(\cdot)$ is odd and the potential
$q(\cdot)$ is even. If $A$ is similar to a self-adjoint operator
then
\begin{equation}
\sup_{\varepsilon>0}\ \frac{\im M_+(i\varepsilon)}{\Real
M_+(i\varepsilon)}<\infty.
\end{equation}
\end{corollary}
\begin{proof}[Proof of Corollary \ref{cor_III_1}]
Since the functions $|r(\cdot)|$ and $q(\cdot)$ are even, one can
easily show that $ m_-(\lambda)=m_+(\lambda)$. It follows from
\eqref{e def psi} that $M_-(\lambda)=-M_+(-\lambda)$. Moreover,
\[
M_+(i\ep)-M_-(i\ep)=M_+(i\ep)+
M_+(-i\ep)=M_+(i\ep)+\overline{M_+(i\ep)} =2\Real M_+(i\ep), \quad
\varepsilon>0.
\]
Combining the last equality with Theorem \ref{necessity}, we
complete the proof.
\end{proof}

\begin{remark}\label{rIII.1}
The case $r(x)= \sgn x$ was studied in papers
\cite{Kar_Mal,KM06}, where several necessary and sufficient
conditions of similarity to a self-adjoint operator have been
obtained. In particular, it was shown that: \emph{if $A$ is similar
to a self-adjoint operator, then
\begin{equation} \label{e kmness}
\frac{\im M_\pm (\eta+i0)}{M_+(\eta+i0)-M_-(\eta+i0)} \in L^\infty
(\R) ,
\end{equation}
where $M_\pm (\eta+i0):= \lim_{\varepsilon \rightarrow +0} M_\pm
(\eta+i\varepsilon)$}.

It is easy to see that condition \eqref{e kmness} is a restriction
of Theorem \ref{necessity} to the real line. Moreover, one can
verify that condition \eqref{e kmness} as well as other necessary
conditions (\cite[Corollaries 5.4--5.6]{KM06}) are fulfilled for
differential operators given in Sections \ref{s ex}--\ref{s exSL},
but the functions \eqref{302} are unbounded in any neighborhood of
zero. Note also that condition \eqref{e kmness} obviously holds if
the operator $A$ has a discrete spectrum. Actually, in this case
$\im M_\pm (\eta+i0) = 0$ a.e. on $\R$.
\end{remark}
\begin{remark}
Let $r(x)=\sgn x$ and $q(\cdot)$ be a finite-zone potential (see
\cite{Lev84}). It was shown in \cite{Kar_Mal,KM06} that the part of the operator $A$ that
corresponds to the essential spectrum of $A$ is similar to a
self-adjoint operator exactly when \eqref{e kmness} holds. Also,
the Jordan structure of the part of $A$ corresponding to the
discrete spectrum has been described in \cite{KM06}. Combining
these results with Theorem \ref{necessity}, we obtain the
following similarity criterion: \emph{the operator $A=(\sgn
x)(-d^2/dx^2 +q(x))$ with a finite-zone potential $q(\cdot)$ is
similar to a self-adjoint operator if and only if the functions
\eqref{302} are bounded on $\C\setminus\R$.}
\end{remark}

%The proofs of Theorem \ref{necessity} and Proposition \ref{p
%NRealSp} are given in Section \ref{s proof}.

\subsection{Applications to J-nonnegative
operators.}\label{ness-positive}

In this subsection we assume that the operator $A$ is $J$-nonnegative.

\begin{proposition} \label{p RealSp}
Let $A$ be the operator associated with \eqref{1}. If $A$ is
J-nonnegative, then the spectrum of $A$ is real.
\end{proposition}
\begin{proof}
Since $L \geq 0$, Proposition \ref{p NRealSp} implies that $\rho
(A) \neq \emptyset $. Theorem \ref{thII.4.1} (i) completes the
proof.
\end{proof}
Theorem \ref{thII.4.1} implies that the operator $A$ admits a spectral function $E (\Delta)$ with
the properties $\textbf{(E1)}-\textbf{(E5)}$. Let us formulate
necessary conditions for regularity of critical points of the
operator $A$.
\begin{theorem}\label{thIV.1}
Let the operator $A$ associated with \eqref{1} be J-nonnegative.
Then:
\begin{description}
\item[(i)] If \ $\infty$ is a regular critical point of $A$, then
for all $R>0$ the functions \eqref{302} are bounded on the set $\{
\lambda \in \C_+ \, : \, |\lambda| > R \}$. \item[(ii)] If \ $0$
is a regular critical point of $A$ and \ $\ker A =\ker A^2$, then
for all $R>0$ the functions \eqref{302} are bounded on the set $\{
\lambda \in \C_+ \, : \, |\lambda| < R \}$.
\end{description}
\end{theorem}
\begin{remark}
It is not hard to see that if $A$ is a J-nonnegative operator and
$\lambda_0 \in \overline{\C_+}\setminus \{0, \infty\}$, then the
functions \eqref{302} are bounded in a sufficiently small
neighborhood of $\lambda_0$.
\end{remark}

\begin{corollary}\label{cIV.1}
Let the operator $A$ associated with \eqref{1} be J-nonnegative.
Assume also that the weight $r(\cdot)$ is odd, the potential
$q(\cdot)$ is even, and $\ker A = \ker A^2$.
If the critical point  $0$ \ ($\infty$) is regular, then %there
%exists a constant $C>0$ such that for all $R>0$ the following
%inequality is valid
\begin{equation}
\im M_+(i\varepsilon)=O(\Real M_+(i\varepsilon)),\quad
\varepsilon\to+0 \qquad (\varepsilon\to+\infty).
\end{equation}
\end{corollary}
The proof is analogous to the proof of Corollary
\ref{cor_III_1} and follows from Theorem \ref{thIV.1}.
%We prove Theorem \ref{thIV.1} and Corollary \ref{cIV.1} in Section
%\ref{s proof}.

In Sections \ref{s ex}-\ref{s exSL}, using the inverse spectral theory of
Sturm-Liouville operators, we will construct the spectral problems
of the form \eqref{1} such that the associated operator $A$ is
J-nonnegative with the singular critical point $0$.

\section{Resolvent estimates \label{s proof}}

\subsection{Boundary triplets for symmetric operators.}

To calculate the spectrum and the resolvent of the operator $A$ we
will use the concepts of boundary triplets %(see \cite{Koch75,Br76,GG})
and abstract Weyl functions (see \cite{DM87,DM95}). Let us briefly
recall basic notions and facts.% concerning this concept.

Let $\H$ and $\HH$ be separable Hilbert spaces. Let $S$ be a
closed symmetric operator in $\H$ with equal and finite deficiency
indices \ $n_+ (S) = n_- (S)=n < \infty$ \ (by definition,
$n_{\pm}(S):=\dim \mathfrak{N}_{\pm i}(S)$, where
$\mathfrak{N}_\lambda(S):=\ker(S^* -\lambda I)$).

\begin{definition}[\cite{GG}]\label{dI.3}
A triplet $\Pi = \{\HH, \Gamma_0,\Gamma_1\}$ consisting of an
auxiliary Hilbert space $\HH$ and linear mappings $\Gamma_j:
\dom(S^*) \longrightarrow \HH$, \ $(j=0,\ 1)$, \ is called a
\emph{boundary triplet for} $S^*$ if the following two conditions
are satisfied:
\begin{description}
\item[(i)] \ $(S^*f,g)_\H-(f,S^*g)_\H = (\Gamma _1 f,\Gamma
_0g)_\HH -        (\Gamma_0 f,\Gamma _1 g)_\HH$, \qquad  $f,\ g
\in \dom(S^*)$; \item[(ii)] the linear mapping $ \Gamma =
\{\Gamma_0 f,\Gamma_1 f\} \ : \dom(S^*) \longrightarrow \HH \oplus
\HH $ is surjective.
\end{description}
\end{definition}
The mappings $\Gamma_0$ and $\Gamma_1$ naturally induce two
extensions $S_0 $ and $S_1$ of $S$ given
by
\[ S_j := S^* \, \upharpoonright \, \dom (S_j), \qquad \dom (S_j)
= \ker \Gamma_j, \quad (j=0,\ 1).
\]
It turns out that $S_0 $ and $S_1$ are self-adjoint operators in
$\H$,\ \ $S_j^*=S_j, \ (j=0,\ 1)$.

The \emph{$\gamma$-field} of the operator $S$ corresponding to the
boundary triplet $\Pi$ is the operator function $\gamma (\cdot) :
\rho(S_0) \to [\HH, \mathfrak{N}_\lambda(S)]$ defined
by $\gamma (\lambda) := (\Gamma_0 \upharpoonright
\mathfrak{N}_\lambda(S))^{-1}$. The function $\gamma $ is
well-defined and holomorphic on $\rho (S_0)$.
\begin{definition}[\cite{DM87,DM95}]\label{dI.6}
Let $\Pi=\{\HH, \Gamma_0,\Gamma_1\}$ be a boundary triplet for the
operator $S^*$. The operator-valued function $M(\cdot) : \rho
(S_0) \to [\HH]$ defined by
$$
M(\lambda) := \Gamma_1 \gamma(\lambda),\qquad \lambda\in\rho
(S_0),
$$
is called \emph{the Weyl function} of $S$ corresponding to the
boundary triplet $\Pi$.
\end{definition}
Note that the Weyl function $M$ is holomorphic on $\rho (S_0)$. It
is well known (see \cite{DM87,DM95}) that the above implicit
definition of the Weyl function is correct and $M(\cdot)$ is an
$(R)$-function obeying $0\in \rho(\im(M(i)))$.

%The Weyl function plays an important role in the spectral theory
%of proper extensions.
%More exactly,
Let $C,D\in[\HH]$. Consider the following extension
$\widetilde{S}$ of $S$, \ \ $S\subset\widetilde{S}$,
\begin{gather}
\widetilde{S}=S_{C,D}:=S^*\upharpoonright\dom(S_{C,D}),\notag \\
\qquad \dom(S_{C,D})=\{f\in\dom(S^*):\
C\Gamma_1f+D\Gamma_0f=0\}.\label{ext_S}
\end{gather}
Notice that each \emph{proper extension} $\widetilde{S}$ of $S$
has the form \eqref{ext_S}, i.e., if $S\subset\widetilde{S}\subset
S^*$, then there exist $C,D\in[\HH]$ such that $\widetilde{S}=S_{C,D}$.

A connection between the Krein—Najmark formula (see, for example,
\cite{AG}) and boundary triplets has been established in
\cite{DM87,DM95}. We use the corresponding result in the form
given in \cite{MalMog}.
\begin{proposition}[\cite{DM87,DM95,MalMog}]\label{p SpKrF}
Suppose $\Pi=\{\mathcal{H}, \Gamma_0, \Gamma_1\}$ is a boundary
triplet for $S^*$, \ $M(\cdot)$ is the corresponding Weyl
function, and $\widetilde{S}=S_{C,D}$, where $S_{C,D}$ is defined
by \eqref{ext_S}. Assume also that $C,\ D,\
(CC^*+DD^*)^{-1}\in[\HH]$. Then:
\begin{description}
\item[(i)] $\lambda \in \rho (S_0)\cap \rho(\widetilde
S) $ \  if and only if \ $ 0 \in \rho (D + C M(\lambda))$.

\item[(ii)] For each $\lambda \in \rho (\widetilde S) \cap \rho
(S_0)$ the following equality holds true
\begin{equation}\label{e KrRes}
(\widetilde S-\lambda)^{-1}
=(S_0-\lambda)^{-1}-\gamma(\lambda)(D +C
M(\lambda))^{-1} C \gamma^*(\overline{\lambda}),
\end{equation}
where the operator-function $\gamma^* (\cdot):
\rho(S_0) \to [\H,\HH]$ is defined by
$$
(\gamma^* (\lambda) f , h)_\HH = (f, \gamma (\lambda) h)_\H,
\qquad f\in\H, \ h \in \HH.
$$
\end{description}
\end{proposition}

\begin{remark}
Assertions (i) and (ii) of Proposition \ref{p SpKrF} are
particular cases of \cite[Corollary 5.3]{MalMog} and
\cite[Corollary 5.6]{MalMog}, respectively. Note that combining
\cite[Proposition 1.6]{DM95}, \cite[Theorem 3.1(1)]{DM95} with
\cite[Lemma 2.1]{MalMog},  one immediately obtains the proof of Proposition
\ref{p SpKrF}.
\end{remark}

\subsection{Boundary triplets for Sturm-Liouville operators.}\label{ss BT}

\textbf{1.} Let $\Aminp$ and $\Aminm$ be the operators
defined in Subsection \ref{ss DifOp}. Since equation \eqref{1} is
in the limit point case at $+\infty$ and $-\infty$, then the
deficiency indices of the symmetric operator $\Aminpm (= \mp
\Lminpm)$ are (1,1) and for all $f,g \in \dom \left(
\Aminpm^* \right)$ we have
\begin{gather}
\left( \Aminpm^* f,\ g \right) -\left( f,\ \Aminpm^*
g \right)= f'(\pm0) \overline{g(\pm0)} -
f(\pm0)\overline{g'(\pm0)} .
\end{gather}
Hence the triplets  $\Pi^+=\{\C, \Gamma_0^+, \Gamma_1^+\}$ and
$\Pi^-=\{\C, \Gamma_0^-, \Gamma_1^-\}$, where
\begin{equation*}%\label{}
\Gamma_0^\pm f:=f'(\pm 0), \quad \Gamma_1^\pm f := -f(\pm 0),
\qquad
% \text{for} \quad
f\in \dom(A_{\min\pm}^*),%\mathfrak{D}_0^{*\pm},
\end{equation*}
are the boundary triplets for $\Aminp^* $ and
$\Aminm^* $, respectively. By the definition of the
functions $\psi_+(\cdot , \lambda)$ and $\psi_-(\cdot , \lambda)$
(see Subsection \ref{sII_2}), we obtain
\begin{equation}\label{308}
\mathfrak{N}_\lambda(\Aminpm)=\ker \left(
\Aminpm^*-\lambda \right) = \{ c \psi_\pm(\cdot , \lambda)
: \ c \in \C \}, \qquad \lambda \in \C\setminus\R.
\end{equation}
Denote by $\gamma^+$ and $\gamma^-$ the $\gamma$-fields
corresponding to the boundary triplets $\Pi^+$ and $\Pi^-$. By
\eqref{e def psi} and \eqref{308}, we get
\begin{gather} \label{e G+def}
\gamma^\pm (\lambda) \, c = \left( \Gamma_0^\pm \, \upharpoonright
\, \mathfrak{N}_\lambda(\Aminpm)
 \right)^{-1}        c
=  c\cdot \psi_\pm(x, \lambda), \qquad c\in\C, \quad \lambda \in
\C\setminus\R.
\end{gather}
Further, the self-adjoint extension
$\Aminpm^*\upharpoonright
\ker(\Gamma_0^\pm)$ of $\Aminpm$ coincides with the operator
$\A0pm$ defined by \eqref{A0_pm} (see Subsection \ref{sII_2}).
The Weyl function $\widetilde M_\pm (\cdot) $  of
$\Aminpm$ corresponding to the boundary triplet
$\Pi^\pm$ is defined by
$$
\widetilde M_\pm (\lambda) := \Gamma_1^\pm \, \gamma^\pm (\lambda)
,\qquad \lambda \in  \rho(\A0pm).
$$
Combining \eqref{e G+def} with \eqref{e def psi}, one obtains
\begin{gather*}
\widetilde M_\pm (\lambda)\, c = \Gamma_1^\pm \, \gamma^\pm
(\lambda) \, c = \Gamma_1^\pm (\, c\, \psi_\pm (\pm0, \lambda) )=
- c \, \psi_\pm (\pm0, \lambda) = c\, M_\pm (\lambda), \qquad
c\in\C, \quad\lambda \in \C\setminus\R .
\end{gather*}
Note that, by definition \ref{dI.6}, the function $\widetilde
M_\pm $ is holomorphic on $\rho (\A0pm)$.  Thus
$\widetilde M_\pm (\cdot) $ is a holomorphic continuation of
$M_\pm (\cdot)$ to the domain $\rho(\A0pm)$. In the
sequel we will write $M_\pm$ instead of $\widetilde{M}_\pm$.

\textbf{2.} Let us consider the symmetric operator $A_{\min}$
defined by \eqref{e Amin=}. Let us determine the linear mappings $
\Gamma_j : \dom(A_{\min}^*)\rightarrow \C^2, \ (j=0,\ 1),$ as
follows
\begin{equation}\label{305}
\Gamma_0 f=\left( \begin{array}{c}
                  f'(+0)         \\
                  f'(-0)
       \end{array} \right) , \qquad
\Gamma_1 f=\left( \begin{array}{c}
                      -f(+0)        \\
                      -f(-0)
                 \end{array} \right) \ ,
\qquad  f\in\dom(A^*_{\min}).
\end{equation}
Since $A_{\min}=\Aminp\oplus \Aminm$ and
$\Gamma_j=\Gamma_j^+\oplus\Gamma_j^-$, then the triplet
$\Pi=\{\C^2, \Gamma_0, \Gamma_1\}$ is a boundary triplet for
$A_{\min}^*$.

Further, we put
\begin{equation}\label{a0}
A_0:=A_{\min}^*\upharpoonright
\ker(\Gamma_0)= \Azp \oplus \Azm.
\end{equation}
Therefore, the operator function $\gamma(\cdot) :
\rho(A_0 ) \to [\C^2, \mathfrak{N}_\lambda(
A_{\min})]$ defined by
\begin{gather} \label{315}
\gamma(\lambda)\left(
\begin{array}{c}
             c_+ \\
             c_-
\end{array} \right):= \gamma^+ (\lambda) c_+ + \gamma^- (\lambda) c_- =
 c_+ \psi_+(\cdot, \lambda) +
 c_- \psi_-(\cdot, \lambda) ,\qquad c_\pm\in\C,
\end{gather}
is the gamma-field corresponding to the boundary triplet $\Pi$.
Moreover, the corresponding Weyl function has the following form
\begin{gather*}
 M(\lambda) := \left(                    \begin{array}{cc}
M_+(\lambda)  &   0  \\
0             &  M_-(\lambda)
\end{array} \right) \ ,
       \qquad
              \lambda \in \rho (A_0) \ \left( = \rho (\Azp ) \cap \rho (\Azm ) \right) .
\end{gather*}

%\subsection{The spectrum and the resolvent of the operator $A$}

%Let $A$, $L$, and $J$ be the operators
%defined in Subsection \ref{ss DifOp}.

\begin{lemma}[cf. \cite{KM06}]
\label{lIII.3} Let $A$ be the operator associated with equation
\eqref{1}, let the operator $A_0$ be defined by
\eqref{a0}. Then:
\begin{description}
\item[(i)] $\sigma(A) \cap \rho(A_0) = \{ \lambda \in
\rho(A_0) \, : \, M_+(\lambda) = M_-(\lambda)\}$;
\item[(ii)] If $\lambda \in \rho(A) \cap \rho(A_0) $,
then the following equality is valid for all $f\in L^2(\R, |r|dx)$
\begin{equation}\label{310}
(A-\lambda )^{-1}f =(A_0 - \lambda )^{-1} f +
\frac{\mathcal{F}_+(f,\lambda)-\mathcal{F}_-(f, \lambda)}
{M_+(\lambda)-M_-(\lambda)} \, (\psi_+(\cdot, \lambda) +
\psi_-(\cdot, \lambda) ),
\end{equation}
where
\begin{gather*}%\label{311}
\mathcal{F}_+(f, \lambda) := \int_0^{+\infty}f(x)\psi_+(x,
\lambda)|r(x)|dx,  \qquad \mathcal{F}_- (f, \lambda):=
\int_{-\infty}^0 f(x) \psi_-(x, \lambda) |r(x)|dx.
\end{gather*}
\end{description}
\end{lemma}

\begin{proof}
\textbf{(i)} Let us rewrite \eqref{e domA} as follows
\begin{gather*}%\label{313}
\dom(A)=\{f\in\dom(A_{\min}^*):
 C \Gamma_1 f + D \Gamma_0 f = 0
\},\\\quad \text{where} \qquad \ C=\left(\begin{matrix}
1&-1\\
0&0
\end{matrix}\right), \ D=\left(\begin{matrix}
                          0&0\\
                          1&-1
                  \end{matrix}\right).
\end{gather*}
By Proposition \ref{p SpKrF} (i), $\lambda\in\rho(A)\cap \rho
(A_0)$ if and only if \ $0\in\rho( D+C M(\lambda))$.
Since
\begin{equation*}
\det(D+C M(\lambda))= \det \left(\begin{array}{cc}
M_+(\lambda)&-M_-(\lambda)\\
1       &      -1
\end{array} \right)=
-M_+(\lambda)+M_-(\lambda),
\end{equation*}
we see that $\lambda\in\rho(A)\cap \rho (A_0)$ exactly
when $M_+ (\lambda) \neq M_- (\lambda)$.

\textbf{(ii)} \
%We use \eqref{e KrRes} to prove assertion (ii).
Let $\lambda \in \rho(A)\cap\rho(A_0)$ and $f \in \H$.
%Note that $\overline{\psi_\lambda^{\, \pm} (\cdot )} = \psi_{\overline{\lambda}}^{\, \pm} ( \cdot)$.
Then after simple calculations we obtain
\begin{gather*}
\gamma^* (\overline{\lambda}) f =
\begin{pmatrix}
\int_0^{+\infty} f(x) \psi_+ (x, \lambda) |r(x)| \, dx\\
\int_{-\infty}^0 f(x) \psi_- (x, \lambda) |r(x)| \, dx
\end{pmatrix}=\begin{pmatrix}
\mathcal{F}_+(f, \lambda)\\
\mathcal{F}_-(f, \lambda)
\end{pmatrix},
\\ \text{and} \qquad
(D+C M(\lambda))^{-1}C = - \frac 1{M_+ (\lambda) - M_- (\lambda)}
\begin{pmatrix}
-1 & 1 \\
-1 & 1
\end{pmatrix}.
\end{gather*}
Combining this equalities with \eqref{315} and \eqref{e KrRes}, we
obtain \eqref{310}.
\end{proof}

\begin{remark}
Lemma \ref{lIII.3} was obtained in \cite{KM06} (for the case $r(x)
= \sgn x$). In the proof given above we use other technique.
\end{remark}

\begin{remark}
Notice also that each point $\lambda_0\in\sigma(A) \cap \rho(A_0)
$ is an eigenvalue of $A$ (see, for example, \cite{AG}), i.e.,
\quad $\sigma(A) \cap \rho(A_0) = \{ \lambda \in \rho(A_0) \, : \,
M_+(\lambda) = M_-(\lambda)\}\subset\sigma_p (A)$.
\end{remark}

\subsection{Proofs.%Regularity of critical points of J-nonnegative operators
\label{ss proofCrP}}

\begin{proof}[Proof of Proposition \ref{p NRealSp}]
Statement \textbf{(i)} obviously follows from Lemma \ref{lIII.3}
and the fact that \\ $\rho (A_0) \subset \R$.

Let us prove \textbf{(ii)}. Assume that the operator $L$ is
semibounded, i.e., $L \geq \eta_0 I, \ \eta_0 \in \R$. Hence
$\Lminp \oplus \Lminm \geq \eta_0 I$ and $\Lminpm \geq \eta_0 I$.
Since $\Aminpm = \pm \Lminpm $, we obtain $\Aminp \geq \eta_0 I$
and $\Aminm \leq -\eta_0 I$.

The operators $\Azp$ and $\Azm$ are self-adjoint extensions of
$\Aminpm$. Furthermore, the deficiency indices of $\Aminpm $ are
$(1, 1)$, hence (see \cite[Chapter VII]{AG}) the operators $\Azp$
and $\Azm$ are semibounded. Therefore, there exists $\eta_1\in
(-\infty, \eta_0]$ such that $\sigma(\Azp)\in [\eta_1, +\infty)$
and $\sigma(\Azm)\in (-\infty, -\eta_1]$. On the other hand, the
operators $\A0pm$ are unbounded. These facts imply $\sigma(\Azp)
\neq \sigma(\Azm)$.

Since $\sigma(\A0pm)=\supp d\tau_\pm$ (see Section \ref{sII_2}),
one immediately gets $\supp d\tau_+ \neq \supp d\tau_-$. By the
Stieltjes inversion formula \eqref{e St}, we conclude that $M_+
(\lambda) \not \equiv M_- (\lambda)$ on $\C\setminus\R$. Hence
Lemma \ref{lIII.3}(i) yields
$$\rho(A) \setminus \R = \{ \lambda \in
\C\setminus\R \, : \, M_+ (\lambda) \neq M_- (\lambda) \} \neq
\emptyset.
$$
\end{proof}

%\subsection{Proof of Theorem
%\ref{necessity}}

The following result is well known.
\begin{proposition}\label{lIII.1}
Let $T$ be a closed operator in a Hilbert space $\mathfrak{H}$ and
$\rho(T)\subset \R$. If $T$ is similar to a self-adjoint operator,
then there exists a positive constant $C>0$ such that
\begin{equation}\label{303}
|\im\lambda|\cdot \|(T-\lambda )^{-1}\|_{\mathfrak{H}}\leq C
\qquad \text{for \ all} \quad \lambda\in\C\setminus\R.
\end{equation}
\end{proposition}
Now we are ready to prove Theorem \ref{necessity}.
\begin{proof}[Proof of Theorem
\ref{necessity}] Suppose that $A$ is similar to a self-adjoint
operator. Then $\sigma(A)\subset \R$. By Lemma \ref{lIII.3}(i), \
$ M_+(\lambda)\neq M_-(\lambda)$ \ for all \ $\lambda \in
\C\setminus \R$. Hence the functions \eqref{302} are well-defined.

Further, by Proposition \ref{lIII.1}, there exists a positive
constant $C>0$ such that
\begin{equation}\label{e RA<C}
|\im\lambda|\cdot \|(A-\lambda )^{-1}\| \leq C \qquad \text{for
all} \quad \lambda\in\C\setminus\R.
\end{equation}
Since the operator $A_0 = A_0^*$ is
self-adjoint, then
$$
|\im \lambda|\cdot\| (A_0 - \lambda)^{-1} \| \leq 1,
\qquad \lambda \in \C\setminus\R.
$$
Combining this inequality with \eqref{e RA<C}, we get
\begin{equation}\label{312}
|\im \lambda|\cdot\left\|(A-\lambda
)^{-1}-(A_0-\lambda)^{-1}
\right\| \leq C+1, \qquad %\text{for} \quad
\lambda\in\C\setminus\R.
\end{equation}
Substituting $f(\cdot)={\psi_\pm(\cdot, \overline{\lambda})}$ in
\eqref{310}, we obtain from \eqref{312} the following inequality
\begin{equation*}
|\im \lambda| \, \frac{\|\psi_\pm(x, \lambda) \| \, (\|\psi_+(x,
\lambda) \|+\|\psi_-(x, \lambda)\|)} {|M_+(\lambda)-M_-(\lambda)|}
\leq \sqrt{2}(C+1), \qquad %\text{for} \quad
\lambda\in\C\setminus \R.
\end{equation*}
Therefore, using \eqref{217}, one immediately gets
\begin{equation*}
\frac{\sqrt{|\im M_{\pm}(\lambda)|} \ (\sqrt{|\im
M_+(\lambda)|}+\sqrt{|\im
M_-(\lambda)|})}{|M_+(\lambda)-M_-(\lambda)|} \leq
\sqrt{2}(C+1) \qquad,
\lambda\in\C\setminus\R.
\end{equation*}
Thus, for \ $ \lambda\in\C\setminus\R$, we have
\begin{equation*}
\frac{|\im M_{\pm}(\lambda)|}{|M_+(\lambda )-M_-(\lambda)|} \leq
\sqrt{2}(C+1) \ .
\end{equation*}
This concludes the proof of Theorem \ref{necessity}.
\end{proof}

\begin{proof}[Proof of Theorem \ref{thIV.1}]

Let us prove assertion $\textbf{(ii)}$. The proof of assertion
$\textbf{(i)}$ is analogous. It is assumed that the operator $A$
has the following properties:
\begin{description}
\item[(A1)] $A$ is a J-self-adjoint J-nonnegative operator;
\item[(A2)]$\ker A =\ker A^2$; \item[(A3)]$0$ is not a singular
critical point of $A$.
\end{description}
Notice that to prove \textbf{(ii)} it is sufficient to show that
the resolvent $(A-\lambda)^{-1}$ of the operator $A$ satisfies the
inequality
\begin{equation}\label{A_e}
|\im \lambda|\cdot\|(A-\lambda)^{-1}\| \leq C \quad \text{for
all}\quad \lambda \in \Omega_R:=\{ \lambda \in \C_+ \, : \, |\lambda| < R \}.
\end{equation}
Actually, if the resolvent of the operator $A$ satisfies
\eqref{A_e}, then, arguing as in proof of Theorem \ref{necessity},
we easily obtain \textbf{(ii)}.

By Proposition \ref{p RealSp} and Theorem \ref{thII.4.1}, the
operator $A$ has a spectral function $E_A(\Delta)$.
%Moreover, by
%\textbf{(A3)} and Theorem \ref{thII.4.2}, $E_A(\Delta)$ is
%well-defined and bounded at zero.
Let us consider the following operator $\mathcal{P}_R:=E_A([-2R,
2R])$, $(R>0)$. Note that the operator $\mathcal{P}_R$ is a
bounded $J$-orthogonal projection in $L^2 (\R, |r(x)|dx)$ (see
\textbf{(E1)-(E2)} in Theorem \ref{thII.4.1}). Furthermore, using
properties \textbf{(E4)--(E5)} of the spectral function
$E_A(\Delta)$, we obtain the decomposition
\begin{equation}\label{dec}
A = \AR \dot + \Ainf , \qquad \AR := A \upharpoonright  \H_0,
\quad \Ainf := A \upharpoonright  \H_{\infty},
\end{equation}
$$
\text{where} \qquad L^2 (\R, |r(x)|dx) = \H_0 \dot + \H_{\infty},\qquad \H_0:= \ran
\left( \mathcal{P}_R \right),
%L^2 (\R, |r(x)|dx)
\qquad \text{and} \qquad \H_{\infty}:= \ran \left( I-\mathcal{P}_R
\right).
$$
 Moreover,
\begin{gather}
\sigma (\AR) \subset [-2R,2R], \qquad \sigma(\Ainf) \subset
(-\infty,-2R] \cup [2R,+\infty) \label{e spAinf}.
\end{gather}
It is obvious that the operator $\AR$ satisfies
\textbf{(A1)--(A3)}. On the other hand, it follows from \eqref{e
spAinf} that $\infty $ is not a critical point of $\AR$. Hence, by
Proposition \ref{p cr=sim}, the operator $\AR$ is similar to a
self-adjoint one. Therefore, by Proposition \ref{lIII.1}, we
obtain
\begin{equation}\label{e resAR}
|\im \lambda|\cdot\|(\AR-\lambda)^{-1}\|_{\H_0} \leq C_1
\quad\text{for all}\quad \lambda \in \C\setminus\R.
\end{equation}
Furthermore,  \eqref{e spAinf} implies
\begin{gather} \label{e resAinf}
\|(\Ainf-\lambda)^{-1}\|_{\H_{\infty}}\leq C_2 \quad \text{for
all}\quad\lambda \in \Omega_R.
\end{gather}
Combining \eqref{e resAR} and \eqref{e resAinf} with \eqref{dec},
we obtain \eqref{A_e}. This completes the proof.
\end{proof}

\section{The operator $-\frac{\sgn x}{(3|x|+1)^{-4/3}} \, \frac{d^2}{dx^2}$ \label{s ex}}

The main aim of Subsections \ref{s ex}--\ref{s exSL} is to present
several explicit examples of indefinite Sturm-Liouville operators
of the form \eqref{3} with the singular critical point $0$. We
start with the case $q\equiv0$. It should be noted that this kind
of operators could be treated by the theory of strings with a
nonmonotone mass distribution function (see \cite{F96}).

%\subsection{The Krein string}
\textbf{1.} In the Hilbert space $L^2(\R, (3|x|+1)^{-4/3}dx)$, let
us consider the operator $A$ defined by the differential
expression
\begin{equation}\label{e eqEx1}
a[y] = -\frac{(\sgn x)}{(3|x|+1)^{-4/3}} \, y''
\end{equation}
on the natural domain $\mathfrak{D}$ (for details see Subsection
\ref{ss DifOp}). Notice that in this case $q(x)\equiv 0$ and
$r(x)=(\sgn x)(3|x|+1)^{-4/3}$.

\begin{theorem}\label{thV_1}
The operator $A$ is J-self-adjoint and J-nonnegative. Moreover,
\begin{description}
\item[(i)] the spectrum of $A$ is real, $\sigma(A)
\subset \R$; \item[(ii)]  $0$ is a simple eigenvalue of $A$;
\item[(iii)] $0$ is a singular critical point of $A$; \item[(iv)]
$A$ is not similar to a self-adjoint operator.
\end{description}
\end{theorem}

The proof of Theorem \ref{thV_1} is based on the following lemma,
which will be proved in the next subsection.
\begin{lemma}\label{lemV_1}
The differential equation
\begin{gather} \label{e eqEx2}
- y''(x) =  \lambda (3x+1)^{-4/3} y(x), \qquad x>0,
\end{gather}
is in the limit point case at $+\infty$. Moreover, the function
\begin{equation}\label{V_1}
m(\lambda)=-\frac{1}{\lambda}+\frac{1}{\sqrt{-\lambda}},
\qquad\lambda\notin\R_+,
\end{equation}
is the Weyl-Titchmarsh $m$-coefficient for \eqref{e eqEx2}.
\end{lemma}
\begin{proof}[Proof of Theorem \ref{thV_1}]
\textbf{(i)} By Lemma \ref{lemV_1}, differential expression
\eqref{e eqEx1} is in the limit point case at both $+\infty$ and
$-\infty$. Hence the operator $A$ is $J$-self-adjoint. Evidently,
the operator $A$ is $J$-nonnegative. It follows from Proposition
\ref{p RealSp} that the spectrum of $A$ is real,
$\sigma(A)\subset\R$. 

Let us prove \textbf{(ii)}. Note that in this case $c(x, 0)\equiv 1$ and $\ker A=\Span\{c(x, 0)\}$.

By Lemma \ref{lemV_1}, we obviously obtain
\begin{equation}\label{W_V_1}
M_+(\lambda)=-\frac{1}{\lambda}+\frac{1}{\sqrt{-\lambda}},\qquad
M_-(\lambda)=-\frac{1}{\lambda}-\frac{1}{\sqrt{\lambda}},\qquad\lambda\in\C\setminus\R.
\end{equation}
Combining \eqref{310} with \eqref{217} and \eqref{W_V_1}, we
obtain the following estimate
\begin{multline}
\|(A-\lambda )^{-1}\| \leq \|(A_0-\lambda )^{-1}\|+
2\cdot\frac{\im M_+(\lambda )+\im M_-(\lambda)}{\im
\lambda\cdot|M_+(\lambda)- M_-(\lambda)|} =\\=
\|(A_0-\lambda
)^{-1}\|+2\cdot\frac{\im(-1/\lambda+1/\sqrt{-\lambda})+\im(-1/\lambda-1/\sqrt{\lambda})}
{\im \lambda\cdot|1/\sqrt{-\lambda}+1/\sqrt{\lambda}|}=\\
=\|(A_0-\lambda I)^{-1}\|+4\frac{\im(-1/\lambda)}{\im
\lambda\cdot|1/\sqrt{-\lambda}+1/\sqrt{\lambda}|}
+2\frac{\im(1/\sqrt{-\lambda}+1/\sqrt{\lambda})}{\im
\lambda\cdot|1/\sqrt{-\lambda}+1/\sqrt{\lambda}|}\leq\\
\leq \|(A_0-\lambda
I)^{-1}\|+\frac{2\sqrt{2}}{|\lambda|^{3/2}}+\frac{2}{|\im
\lambda|}, \qquad  \im \lambda\neq 0.
\end{multline}
Hence, for $\ep>0$
$$
\|(A-i\ep )^{-1}\|=O(|\ep |^{-3/2}) \qquad \ep \to 0.
$$
Therefore, the Riesz index of $0$ is less than 2, i.e.,
$\ker A^2=\ker A$.

To prove \textbf{(iii)} we use Corollary \ref{cIV.1}.

Simple calculations show that
\[
\im M_+(i\ep )=\frac{1}{\ep}+\frac{1}{\sqrt{2\ep}}, \qquad \Real
M_+(i\ep)=\frac{1}{\sqrt{2\ep}}, \qquad \ep>0,
\]
and
\begin{equation}
\frac{\im M_+(i\ep)}{\Real
M_+(i\ep)}=\frac{1/\ep+\sqrt{1/2\ep}}{\sqrt{1/2\ep}}=1+\sqrt{\frac{2}{\ep}}\to+\infty,\quad
\ep\to+0.
\end{equation}
Thus the condition of Corollary \ref{cIV.1} fails, hence 0 is a
singular critical point of the operator $A$.

Finally, notice that \textbf{(iv)} directly follows from
\textbf{(iii)} (see Proposition \ref{p cr=sim}).
\end{proof}

%%%%%%%%%%%%%%%%%%%%%%%%%%%%%%%%%%%%%%%%%%%%%%%%%%%%%%%%%%%%
\textbf{2.} Let us briefly recall the basic facts from the Krein
string spectral theory (see \cite{KK2}, and also
\cite{DymMcKean}).

A string $\mathcal{S}$ is specified by a pair $\len$ and
$\m$, where the number $\len>0$ is the length of
$\mathcal{S}$, and the function $\m: [0, \len)\to \R_+$ is
the mass distribution. Naturally, $\m$ is nonnegative,
nondecreasing, continuous from the right, and $x=0$ is a point of
growth. With such a string one can naturally associate a
self-adjoint operator $\L$ (see \cite{DymMcKean}). This is done
by restricting the formal differential operator $\ell:f\to
-d^2f/d\m dx$ to the special domain $\dom(\L)$ in the Hilbert
space $L^2([0,\len], d\m)$. If
\[
\m \in AC_{loc} , \quad
d\m(x)=r(x)dx, \quad \len=+\infty, \quad \text{and} \quad
\int_{\R_+} t^2d\m(t)=\infty,
\]
then $\L$ coincides
with the Sturm-Liouville operator $\Azp$ defined by \eqref{A0_pm}.

It should be noted that for the string $\mathcal{S}$ one can
naturally determine the Weyl-Titchmarsh $m$-function $m(\cdot)$
and the spectral function $\tau(\cdot)$. Note that if $\m$ is
locally absolutely continuous, $d\m(x)=r(x)dx$, then $m$ and
$\tau$ become the classical $m$-function and the classical
spectral function.

The following fundamental result is due to M. Krein \cite{krein52}
(see also \cite[\S 11]{KK2} and \cite[\S 6.6]{DymMcKean}).
\begin{theorem}[\cite{krein52}]\label{thII.2.1}
A nondecreasing function $\tau:\R_+\to\R_+, \ (\tau(0)=0)$ is a
spectral function of a string  if and only if
\begin{equation}\label{225}
\int_0^{+\infty}\frac{d\tau(s)}{1+s}<\infty.
\end{equation}
Under this condition a string $\mathcal{S}$ (i.e., the length $\len$
and the mass distribution $\m$) is uniquely determined by $\tau$.
\end{theorem}
In the following we also need  rule for the change in the string
$\mathcal{S}\to \mathcal{S}^*$ resulting from a change
$\tau\to\tau^*$ in its spectral function. This fact has been
discovered by M. Krein \cite[Theorem 2.3]{krein53} %and goes back
%to Stieltjes \cite[\S 10]{St}
(see also \cite[\S 6.9, Rule 2]{DymMcKean}).
\begin{theorem}[\cite{krein53}]\label{lII.2.1}
Let $\m(x)$ and $\m^*(x)$ be the mass distributions of two
different strings $\mathcal{S}$ and $\mathcal{S}^*$ with lengths
$\len$ and $\len^*$ respectively. Let $\tau$ and $\tau^*$ be the
spectral functions of $\mathcal{S}$ and $\mathcal{S}^*$. If
$c>-\rho(+0)$ and $\tau^*(s)=c+\tau(s)$ for all $s>0$, then
\begin{equation}\label{226}
\m^*(x)=\frac{\m(\zeta)}{1+c\m(\zeta)}, \qquad x=\int_0^\zeta
(1+c\m(s))^2ds, \quad (0\leq \zeta \leq \len).
\end{equation}
\end{theorem}
%%%%%%%%%%%%%%%%%%%%%%%%%%%%%%%%%%%%%%%%%%%%%%%%
Now we are ready to prove Lemma \ref{lemV_1}.
\begin{proof}[Proof of Lemma \ref{lemV_1}]
Putting $\lambda=0$ in \eqref{e eqEx2} , we obtain
\[
c(x, 0) = 1,\qquad s(x, 0)=x,\qquad x>0.
\]
Since $s(x, 0)\notin L^2(\R_+, (3x+1)^{-4/3})$, the Weyl
alternative implies that expression \eqref{e eqEx2} is in the limit point case
at $+\infty$.

Further, one can easily compute that
\[
m(\lambda)=-\frac{1}{\lambda}+\frac{1}{\sqrt{-\lambda}}=\int_0^{+\infty}\frac{d\tau(s)}{s-\lambda},
\qquad \lambda\notin\R_+,
\]
with
\begin{equation}\label{s1}
\tau(s):=\left\{ \begin{array}{cc}
   1+2\sqrt{s}/\pi, & s>0            \\
  0   , &  s\leq 0
                             \end{array} \right.  \ .
\end{equation}
By Theorem \ref{thII.2.1}, $\tau(\cdot)$ is a spectral function of
a certain string $\mathcal{S}$. Let us recover the mass distribution
$\m$ and the length $\len$ of $\mathcal{S}$.

It is well-known that the function
\begin{equation}\label{s0}
\tau_0(s):=\left\{ \begin{array}{cc}
   2\sqrt{s}/\pi, & s>0            \\
  0   , &  s\leq 0
                             \end{array} \right.  \ ,
\end{equation}
is the spectral function of the problem
\[
-y''(x)=\lambda y(x),\qquad x>0;\qquad y'(0)=0.
\]
In other words, $\tau_0(\cdot)$ is a spectral function of the string
$\mathcal{S}_0$ with the mass distribution $\mathcal{M}_0(x)=x$ and the length $l_0=+\infty$. Using
Theorem \ref{lII.2.1}, we obtain
\[
\mathcal{M}(x)=\frac{\mathcal{M}_0(\zeta)}{1+\mathcal{M}_0(\zeta)}=\frac{\zeta}{1+\zeta},
\qquad x%%=\int_0^\zeta (1+\mathcal{M}_0(s))^2ds
=\int_0^\zeta (1+s)^2ds=\frac{(1+\zeta)^3-1}{3}.
\]
Hence,
\begin{equation}
\zeta=(3x+1)^{1/3}-1, \qquad \mathcal{M}(x)=1- (3x+1)^{-1/3},
\quad 0\leq x<l=+\infty.
\end{equation}
Finally, $\mathcal{M}(x) $ is locally absolutely continuous on
$\R_+$ and
\[
d\mathcal{M}(x)=r(x)dx = (3x+1)^{-4/3}dx , \qquad x>0 .
\]
Therefore, the function \eqref{s1} is a spectral function of the
boundary value problem \eqref{e eqEx2}.
\end{proof}

\section{Operators with the singular critical point zero: \\
the case $r(x)=\sgn x$ \label{s exSL}}

In this section we suppose that $r(x) = \sgn x$, $x \in \R$. Two
examples of indefinite Sturm-Liouville operators of type
$(\sgn x)(-d^2/dx^2 +q)$ will be considered.

\subsection{An operator with a decaying potential.}\label{sub_6_1}

The main object of this subsection is the following
%indefinite Sturm-Liouville
operator
\begin{equation}\label{VI_1_01}
(Ay)(x)=(\sgn
x)\left(-y''(x)+6\frac{x^4-6|x|}{(|x|^3+3)^2}y(x)\right),\qquad
\dom(A)=W_2^2(\R).
\end{equation}
Here $W_2^2(\R)$ is the Sobolev space. Notice that the potential
$q(x)=6(x^4-6|x|)(|x|^3+3)^{-2}$ is bounded on $\R$, hence the
operator $A$ is $J$-self-adjoint.

\begin{theorem}\label{thVI_1}
Let $A$ be the operator of the form \eqref{VI_1_01}. Then:
\begin{description}
\item[(i)] $A$ has a real spectrum, $\sigma (A) \subset \R$;
\item[(ii)] $0$ is a simple eigenvalue of $A$; \item[(iii)] $A$ is
not similar to a self-adjoint operator.
\end{description}
\end{theorem}
As in the previous section, we start with a preliminary lemma.
\begin{lemma}\label{lemVI_1}
The function
\begin{equation}\label{VI_I_02}
m_0(\lambda)=\frac{\lambda}{1+\lambda\sqrt{-\lambda}}, \qquad
\lambda\notin\R,
\end{equation}
is the Weyl-Titchmarsh $m$-coefficient for the boundary value problem
\begin{equation}\label{VI_I_03}
-y''(x)+6\frac{x^4-6|x|}{(|x|^3+3)^2}y(x)=\lambda\ y(x),\qquad
x\geq 0;\quad y'(0)=0.
\end{equation}
\end{lemma}
\begin{proof}[Proof of Lemma \ref{lemVI_1}]
Let us consider the following function
\begin{equation}\label{s3}
 \tau_\infty (s):=\left\{ \begin{array}{cc}
                          1+\frac{2}{3\pi}s^{3/2}, & s>0            \\
                            0              ,       & s\leq 0
                             \end{array} \right.  \  .
\end{equation}
Using the algorithm of Gelfand and Levitan (see e.g. \cite{Lev84}), we
obtain that (calculations are omitted) $\tau_\infty$ is the
spectral function of the problem
 \begin{gather}\label{VI_I_04}
-y''(x)+6\frac{x^4-6|x|}{(|x|^3+3)^2}y(x)=\lambda\ y(x),\qquad
x\geq 0;\quad y(0)=0.
\end{gather}

Moreover, the function
\[
m_\infty(\lambda):=-\frac{1}{\sqrt{2}}+\int_{0}^{+\infty}\left(\frac{1}{s-\lambda}-\frac{\lambda}{1+\lambda^2}\right)d\tau_\infty(s)=-\frac{1}{\lambda}-\sqrt{-\lambda},\qquad
\lambda\notin\R_+,
\]
is the Weyl-Titchmarsh $m$-coefficient for the boundary value
problem \eqref{VI_I_04}. But, it is obvious that
$m_0(\lambda)=-1/m_\infty(\lambda)$,\
$\lambda\notin\{-1\}\cup\R_+$. Thus (see Remark \ref{remII_I_1})
$m_0$ is the Weyl-Titchmarsh $m$-coefficient for the problem
\eqref{VI_I_03}.
\end{proof}

\begin{proof}[Proof of Theorem \ref{thVI_1}]
\textbf{(i)} By Lemma \ref{lemVI_1}, we have
\begin{equation} \label{e M+=M-}
M_+(\lambda)=m_0(\lambda)=\frac{\lambda}{1+\lambda\sqrt{-\lambda}},
\qquad M_-(\lambda)=-M_+(-\lambda),\qquad \lambda\notin\R.
\end{equation}
Since
\[
\frac{1}{M_+(\lambda)}-\frac{1}{M_-(\lambda)}=-\frac{1}{\lambda}-\sqrt{-\lambda}-\left(-\frac{1}{\lambda}+\sqrt{\lambda}\right)
=-\sqrt{-\lambda}-\sqrt{\lambda}\neq 0 \quad \text{for all}\quad
\lambda\notin\R,
\]
we see that $M_+(\lambda)\neq M_-(\lambda)$ for
$\lambda\notin\R$.
By Lemma \ref{lIII.3}, the operator $A$ has a
real spectrum.

\textbf{(ii)} Since $\tau_\infty(\cdot)$ is the spectral function
of the boundary value problem \eqref{VI_I_04} and
$\tau_\infty(+0)\neq\tau_\infty(-0)$, then $\lambda=0$ is an
eigenvalue of the problem \eqref{VI_I_04}. This implies $s(x,
0)\chi_+(x)\in L^2(\R_+)$ (for the definition of $s(x, \lambda)$
see Subsection \ref{sII_2}). Moreover, the potential $q$ is even.
Hence, $s(x, 0)\chi_-(x)\in L^2(\R_-)$. Thus $s(x, 0)\in L^2(\R)$
and therefore $s(x, 0)\in \ker L$.  Since $\ker A=\ker
L\neq\{0\}$, we conclude that $0\in \sigma_p(A)$. Notice also that
\eqref{VI_I_04} is limit point at $+\infty$. Therefore, $c(x,
0)\notin L^2(\R)$ and $\ker A=\ker L=\Span\{s(x,0)\}$.

To prove that $\ker A=\ker A^2$ we check the conditions of
Proposition \ref{p RindEig}. It is not so hard to obtain that
\[
M_+(\lambda) =
m_0(\lambda)=-\frac{2}{3(1+\lambda)}+\int_{0}^{+\infty}\frac{d\tau_{ac}(s)}{s-\lambda},\qquad
\lambda\notin \{-1\}\cup\R_+,
\]
with
\[
 d\tau_{ac}(s):=\frac{s^{5/2}}{\pi(1+s^3)}ds,\qquad s>0.
\]
Thus $d\tau_+ (s)= \frac{2}{3} \, \delta (s+1)ds+ d\tau_{ac} (s)$,
where $\delta (s)ds$  is the Dirac measure. It follows from
\eqref{e M+=M-} that $d\tau_- (s)= \frac{2}{3} \, \delta (s-1)ds -
d\tau_{ac} (-s)$. Hence,
\begin{gather*}
\tau_\pm (+0) = \tau_\pm (-0), \qquad \int_{\R } \frac 1{|s |^2} \ d \tau_+ (s)  = \int_{\R } \frac 1{|s
|^2} \ d \tau_- (s) = \frac 2{3 } + \int_{\R } \frac {s^{1/2}}{\pi
(1+s^3)} \ d s < \infty,
\end{gather*}
and
\begin{gather*}
\int_\R \frac {d\tau_+ (s)}{s}  = M_+ (0+i0) = 0 = M_- (0+i0) =
\int_\R \frac {d\tau_- (s)}{s} .
\end{gather*}
So conditions \eqref{e tau=tau} and \eqref{e fr1=fr1} are
fulfilled and, by Proposition \ref{p RindEig}, we have $0 \in \sigma_p
(A)$ and $\dim(\ker A)=1$. Since
\[
\int_\R\frac{d\tau_\pm (s)}{s^4} >
\int_\R\frac{d\tau_{ac}(s)}{s^4}
=%\int_0^{+\infty}\frac{s^{5/2}ds}{\pi(1+s^3)s^4}=
\int_0^{+\infty}\frac{ds}{\pi(1+s^3)s^{3/2}}=\infty \  ,
\]
we see that both conditions \eqref{e frj<inf r} fail. Thus,
Proposition \ref{p RindEig} yields that $0$ is a simple eigenvalue
of the operator $A$.

 \textbf{(iii)} After simple calculations, one gets for $\ep>0$
\[
\im\frac{1}{M_+(i\ep)}=\frac{1}{\ep}+\frac{\sqrt{\ep}}{\sqrt{2}} \ , \qquad
\Real\frac{1}{M_+(i\ep)}=-\frac{\sqrt{\ep}}{\sqrt{2}}.
\]
Hence,
\begin{gather*}
\frac{\im M_+(i\ep)}{\Real M_+(i\ep)}=\frac{\im(1/ M_+(i\ep))}{
\Real(1/M_+(i\ep))}=\frac{1/\ep +\sqrt{\ep/2}}{-\sqrt{\ep/2}}\to
-\infty, \qquad \ep\to +0.
\end{gather*}
Therefore, by Corollary \ref{cor_III_1}, the operator $A$ is not
similar to a self-adjoint one.
\end{proof}
\begin{remark}
Note that the operator $A$ of the form \eqref{VI_1_01} is not
$J$-nonnegative, but it is definitizable. Actually, consider the
corresponding self-adjoint operator $L:=JA$. It is easy to show
that
%Since $m_0$ and
%$m_\infty$ are analytic on $(-\infty, -1)\cup(-1, 0)$ and the
%potential $q(\cdot)$ is even, then
%the spectrum of $L$ is
$\sigma(L)=\{-1\}\cup[0, +\infty)$ and $\lambda_0 = -1$ is an
eigenvalue of $L$. Therefore, the form $[A\cdot, \cdot\ ]=(L\cdot,
\cdot\ )$ has exactly one negative square. Hence the operator $A$
is definitizable (see \cite[Proposition 1.1]{CurLan}).

Using the arguments from the proof of Theorem \ref{thIV.1}, it is
not hard to show that $0$ is a singular critical point of the
operator $A$.
\end{remark}
%%%%%%%%%%%%%%%%%%%%%%%%%%
\subsection{$J$-nonnegative Sturm-Liouville operator with the singular critical point zero.}
\label{sub_6_2} \qquad

The following result follows easily from \cite[Lemma 3.5
(iii)]{CurLan} and \cite[Theorem 3.6 (i)]{CurLan}: if the operator
$L=-d^2/dx^2+q(x)$ (acting in $L^2 (\R)$) is nonnegative, then
$\infty$ is a regular critical point of the operator $A=(\sgn
x)L$.
%In \cite{CurLan}, \'Curgus and Langer shown that $\infty$ is a
%regular critical point of a $J$-nonnegative Sturm-Liouville
%operator $A$ of the form \eqref{3} if $r(x)=\sgn x$.
%In other words, a $J$-nonnegative Sturm-Liouville operator
%$A=(\sgn x)(-d^2/dx^2+q(x))$ with a bounded potential is not similar to a self-adjoint operator if
%and only if $0$ is a singular critical point of $A$ or $\ker
%A\neq\ker A^2$.
The regularity of critical point $0$ of the operator $A=(\sgn x)(-d^2/dx^2+q(x))$
was proved for the following three cases:
\begin{description}
\item[(i)] $q \equiv 0$ (see \cite{CN95});
\item[(ii)] the spectrum of $A$ is real and $q$ satisfies condition \eqref{e fmoment}
(see \cite{FSh1});
\item[(iii)] $A$ is definitizable and $q$ is a finite-zone potential (see \cite{Kar_Mal,KM06}).
\end{description}

The goal of this subsection is to
show that there exists a J-nonnegative operator of the type
$A=(\sgn x)(-d^2/dx^2+q(x))$ with the singular critical point $0$.

We first need in some preparations. Let us consider the following
nondecreasing function
\begin{equation}\label{s2}
\tau(s):=\left\{ \begin{array}{cc}
                          1+\frac{2}{\pi}(\sqrt{s}-\arctan\sqrt{s}), & s>0            \\
                            0 ,             &  s\leq 0
                             \end{array} \right.  \  .
\end{equation}
By the Gelfand--Levitan theorem (see \cite{LevSar70}),
$\tau(\cdot)$ is a spectral function of the
boundary value problem
\begin{equation}\label{VI_II_08}
-y''(x)+q_0(x)y(x)=\lambda y(x), \quad y'(0)=0, \qquad x\in[0,
+\infty),
\end{equation}
with a certain continuous potential $q_0$.
It follows from \eqref{216} that the corresponding Weyl-Titchmarsh
$m$-coefficient has the form
\begin{equation}\label{W2}
m(\lambda)=\int_{-0}^{+\infty}\frac{d\tau(s)}{s-\lambda}=-\frac{1}{\lambda}+\frac{1}{\sqrt{-\lambda}}-\frac{1}{-\lambda+\sqrt{-\lambda}},\qquad
\lambda\notin[0, +\infty).
\end{equation}
Let us recover the corresponding differential expression, i.e.,
the potential $q_0(x)$. Using the Gelfand--Levitan algorithm, we obtain
\begin{equation}\label{II_72_05}
f(x,y)=\int_0^{+\infty}\cos\sqrt{\lambda}x\ \cos\sqrt{\lambda}y\
d(\tau(\lambda)-2\sqrt{\lambda}/\pi)=1-e^{-x}\cosh y, \qquad 0\leq
y\leq x,
\end{equation}
and
\begin{equation}\label{II_72_07}
q_0(x)=2\frac{d}{dx}K(x,x),\qquad x>0,
\end{equation}
where the kernel $K(x, y)$ is the solution of the
Gelfand--Levitan equation
\begin{equation}\label{II_72_06}
K(x,y)+f(x,y)+\int_0^x K(x,t)f(t,y)dt=0,\qquad 0\leq y\leq x.
\end{equation}

\begin{theorem}\label{thII_7_2}
Let the potential $q_0$ be defined by
\eqref{II_72_07}, \eqref{II_72_06}, and \eqref{II_72_05}. Let the operator $A$ be defined by the differential
expression
\begin{equation}\label{II_72_10}
(\sgn x)\left(-\frac{d^2}{dx^2}+q_0(|x|)\right)
\end{equation}
on the natural domain $\mathfrak{D}$ in the Hilbert space
$L^2(\R)$ (for the definition see Subsection \ref{ss DifOp}).
Then:
\begin{description}
\item[(i)] $A$ is a $J$-nonnegative J-self-adjoint operator and
$\sigma(A)\subset \R$;

\item[(ii)]  $0$ is a simple eigenvalue of $A$;

\item[(iii)]  $0$ is a singular critical point of $A$;

\item[(iv)] the operator $A$ is not similar to a self-adjoint
operator.
\end{description}
\end{theorem}

\begin{proof}
By Proposition \ref{p RealSp}, to prove $\textbf{(i)}$ we only have to show
that the operator $L= -d^2/dx^2+q_0(|x|)$ is nonnegative.
Combining arguments from Lemma \ref{lIII.3} (i) with Proposition \ref{p SpKrF} (i),
one can show that
\begin{gather}\label{e sL}
\sigma(L) \cap \rho(\Lzp \oplus \Lzm ) = \{ \lambda \in
\rho(\Lzp \oplus \Lzm ) \, : \, m_+(\lambda) + m_-(\lambda) = 0\} .
\end{gather}
Since the potential $q_0 (|\cdot|)$ is even, we see that
$m_-(\lambda)=m_+(\lambda)=m(\lambda)$. Moreover, \eqref{e sA=supp} implies
that
\begin{gather*} \label{e slpm}
\sigma (\Lzp)= \sigma (\Lzm) = \supp d\tau = [0,+\infty) .
\end{gather*}
From this, \eqref{e sL} and \eqref{W2}, we obtain
$\sigma(L) \subset \R_+$, i.e., $L \geq 0$.
(It is not difficult to show that $\sigma(L) = \R_+$.)

$(ii)$ Since $\tau$ is the spectral function of the problem
\eqref{VI_II_08} and $\tau (+0) \neq \tau (-0)$, we see that
$\lambda=0$ is an eigenvalue of the problem \eqref{VI_II_08}. That
is $c(x, 0)\chi_+(x)\in L^2(\R_+)$. Furthermore, the potential
$q_0(|x|)$ is even, hence $c(x, 0)\chi_-(x)\in L^2(\R_-)$ and
$c(x, 0)\in\ker L$. Let us note that the operator $L$ is
self-adjoint, i.e., the differential equation \eqref{e eq|r|} with
$r(x)=\sgn x$ and $q(x)=q_0(|x|)$ is limit point at $+\infty$ and
$-\infty$. Therefore, $s(x,0)\notin L^2(\R)$ and $\ker
L=\Span\{c(x, 0)\}$.

The equality $\ker A=\ker L$ implies $0\in \sigma_p(A)$ and
$\dim(\ker A)=1$. Arguing as in the proof of Theorem \ref{thV_1}, one can show that $\ker A=\ker A^2$, so
$0$ is a simple eigenvalue of $A$. On the other hand, this fact follows from
\cite[Theorem 1 (2.ii)]{KarKr05} (see also \cite[Theorem 4.2 (2.ii)]{KM06}) since
$
\int_{+0}^{\infty}s^{-2}d\tau(s)=+\infty
$.
This completes the proof of $(ii)$.

 $(iii)$ By \eqref{W2}, we obtain
$$
\im
m_+(iy)=\frac{1}{y}+\frac{1}{\sqrt{2y}}-\frac{y+\sqrt{y/2}}{y/2+(y+\sqrt{y/2})^2}\geq
\frac{1}{y},\qquad y\geq0;
$$
$$
\Real m_+(iy) = \frac{1}{\sqrt{2y}} -
\frac{\sqrt{y/2}}{y/2+(y+\sqrt{y/2})^2} =
\frac{1+\sqrt{y/2}}{1+y+\sqrt{2y}},\qquad y\geq0;
$$
Combining these relations, one easily gets
$$
\lim_{y\to +0}\frac{\im m_+(iy)}{\Real m_+(iy)}=+\infty.
$$
Therefore, by Corollary \ref{cIV.1}, $0$ is a singular critical
point of the operator $A$.

$(iv)$ Follows from Proposition \ref{p cr=sim} and $(iii)$.
\end{proof}

\begin{remark}
It is easy to see that the potential $q_0$ is continuous. 
%Notice that the potential $q_0$ 
%belongs to $C_{loc}^\infty(\R_+)$ since the function
%$$
%\Phi(x)=\int_{-\infty}^{+\infty}\cos(\sqrt{\lambda}x)
%d\tau(\lambda)=1-\cosh x,\qquad x\geq0,
%$$
%is in $C_{loc}^\infty(\R_+)$ (see \cite{krein53b}).
Moreover, the potential $q_0$ is $L^2$ potential, $q_0\in
L^2(\R)$. We plan to publish a proof of this fact in the
forthcoming paper devoted to indefinite Sturm-Liouville
problems with decaying potentials.
\end{remark}

\ack{The authors are indebted to Andreas Fleige for bringing the
problem to their attention. The authors are deeply grateful to
Mark Malamud for his constant attention to this work and numerous
fruitful discussions. The first author would like to thank Michel Chipot for the hospitality 
of the University of Zurich.}

%%%%%%%%%%%%%%%%%%%%%%%%%%%%%%%%%%

%\affiliation{Department of Nonlinear Analysis,
%          Institute of Applied Mathematics and Mechanics,
%          National Academy of Sciences of Ukraine,
%          R. Luxemburg str., 74, Donetsk 83114, UKRAINE
%\email{duzer80@mail.ru}          }
\end{document}